\documentclass[11pt,reqno]{amsart}
\usepackage[dvips]{graphics}
\usepackage{amssymb}
\usepackage{dsfont}

\numberwithin{equation}{section}

\renewcommand{\qed}{\hspace*{\fill}q.e.d.\vspace{2ex}}

\newtheorem{thm}{Theorem}[section]
\newtheorem*{thm*}{Theorem}
\newtheorem*{thmmain*}{MAIN THEOREM}
\newtheorem{lem}[thm]{Lemma}
\newtheorem{cor}[thm]{Corollary}
\newtheorem{prop}[thm]{Proposition}
\newtheorem*{prop*}{Proposition}
 
\theoremstyle{definition}

\theoremstyle{remark}
\newtheorem{rem}{Remark}[section]
\newtheorem{ex}[rem]{Example}

\newcommand{\tref}[1]{Theorem~\ref{#1}}
\newcommand{\cref}[1]{Corollary~\ref{#1}}
\newcommand{\pref}[1]{Proposition~\ref{#1}}

\newcommand{\lref}[1]{Lemma~\ref{#1}}

\def\dim{\mathop{\text{dim}}}
\def\ind{\mathrm{\text{ind}}}
\def\codim{\mathop{\text{codim}}}
\def\Id{\mathop{\text{Id}}}
\def\Jac{\mathrm{\text{Jac}}}

\begin{document}

\title{Curvature explosion in quotients and applications}

\author{Alexander Lytchak}
\address{A. Lytchak, Mathematisches Institut, Universit\"at Bonn,
Beringstr. 1, 53115 Bonn, Germany}
\email{lytchak\@@math.uni-bonn.de}

\author{Gudlaugur Thorbergsson}
\address{G. Thorbergsson, Mathematisches Institut,
Universit\"at zu K\"oln,
Weyertal 86 - 90, 
50931 K\"oln, Germany }
\email{gthorbergsson\@@mi.uni-koeln.de}



\begin{abstract}
We prove that the quotient space of a variationally complete group action
is a good Riemannian orbifold. The result is generalized to
singular Riemannian foliations without horizontal conjugate points. 
\end{abstract}

\thanks{The first author was supported in part by the SFB 611 
{\it Singul\"are Ph\"anomene und Skalierung in mathematischen Modellen.} The
second author was supported in part by the DFG-Schwerpunkt 
{\it Globale Differentialgeometrie.}}

\maketitle
\renewcommand{\theequation}{\arabic{section}.\arabic{equation}}
\pagenumbering{arabic}

\section{Introduction}

Let $M$ be a Riemannian manifold and let $K$ be a closed group of isometries
of $M$. Usually, the quotient space $B=M/K$ is not a Riemannian manifold,
but an Alexandrov space of curvature (locally) bounded below,
 stratified by Riemannian  manifolds  that consist of orbits of the same type.
The set $M_0$ of all points in $M$ with principal isotropy group
is open and dense in $M$, and it is invariant under the action of
 $K$. The quotient $B_0=M_0 /K$
is the maximal stratum of $B$. In many cases,  the  
sectional curvatures in $B_0$ explode as one approaches a singular point  
$y\in B_{sing} =B\setminus B_0$. However, sometimes it does not happen, 
like in the  case  of an exceptional orbit. An  interesting class
of examples is given by polar actions, where the quotient is a smooth 
Riemannian
orbifold and, therefore, its sectional curvatures are uniformly bounded
on compact subsets. Our first objective is a precise description
of such (non-) explosions.
For a point $z\in B_0$, we denote by $\bar \kappa (z)$ the maximum
of the sectional curvatures of tangent planes at $z$. Then we have:

\begin{thm} \label{theorem1}
Let $M$ be a Riemannian manifold and let $K$ be a closed group of isometries
of $M$. Let $B=M/K$ be the quotient. Let $x\in M$ be a point with isotropy
group $K_x$ acting on the normal space $H_x$ of the orbit $Kx \subset M$.
Set $y=Kx \in B$. Then the following are equivalent:

\begin{enumerate}
 \item $\limsup _{z\in B_0, z\to y}  \bar \kappa (z) < \infty $;

 \item  $\limsup _{z\in B_0, z\to y}  \bar \kappa (z) \cdot  d^2 (y,z) =0$;

 \item The action of $K_x$ on $H_x$ is polar;

 \item A neighborhood of $y$ in $B$ is a smooth Riemannian orbifold.  
\end{enumerate}
 \end{thm} 

\begin{rem}
With arguments similar to those used in the proof of
\tref{theorem1} one sees  that the explosion of curvatures in $B_0$
is at most quadratic in the distance to the singular strata; namely, for each
compact subset $U$ of $B$ there is a constant $C>0$ such that 
for all $z\in B_0 \cap U$
$$
\bar \kappa (z) \leq \frac{C}{ d^2  (z, B_{sing})}.
$$
\end{rem}

Since all orthogonal actions on Euclidean spaces of cohomogeneity at most
$2$ are polar, we deduce that the complement of the  
set of Riemannian orbifold points in the quotient  has codimension at least
$3$:  

\begin{cor} \label{strata}
 Let $M$ be a Riemannian manifold and let $K$ be a closed group of isometries
of $M$.  Then each point $y$ in the quotient $B=M/K$   contained in a stratum
of codimension $\leq 2$ has  a neighborhood in $B$ isometric to a smooth 
Riemannian orbifold. 
\end{cor}

 This result makes large parts of the quotient accessible with differential 
geometric tools and can be used to derive some geometric
properties of the quotient, like the invariance of the Liouville measure under
the (quasi-) geodesic flow; see \tref{quasi} below.

  \tref{theorem1} can be applied  to study variationally complete actions.
Variationally complete 
 actions were introduced  by Bott and used in \cite{Bott} and
 \cite{BottS}, as a mean
for studying  the topology of (loop spaces of) symmetric spaces. 
In these papers it is  observed that the 
orbits of such actions are taut submanifolds of the ambient (symmetric)
space, which establishes  strong relations
between the topology of the ambient manifold and the topology of the orbits.
 Conlon (\cite{Conlon}) showed  that hyperpolar actions
 are variationally complete.  The converse
was proven in \cite{DO} for Euclidean spaces, in \cite {GT} for compact
symmetric spaces
and, more generally, in \cite{LT} for 
non-negatively curved Riemannian manifolds.

  In general, a variationally complete action of a group
$K$ on a Riemannian manifold $M$ does not need to be  
 polar, nor does an  action with taut orbits need
to be variationally complete.   We describe all variationally complete
actions in terms of the quotient space and give a precise meaning to 
the following remark of Bott and Samelson
(\cite{BottS}, p.965):  ``Intuitively, we like to think of
variational completeness as absence of conjugate points on the decomposition
space $M/K$.''

\begin{thm} \label{mainthm}
Let $M$ be a complete Riemannian manifold and $K$ a closed group
of isometries of $M$. The action of $K$ on $M$ is variationally complete
if and only if the quotient $B=M/K$ is a good Riemannian 
 orbifold without conjugate points, i.e., if and only if $B$ is isometric to a 
quotient $B=N/\Gamma$, where $N$ is a smooth, complete, simply connected
Riemannian manifold  without conjugate points and $\Gamma$ a discrete group
of isometries of $N$.
\end{thm}

 If the manifold $M$ is non-negatively curved, then so is the quotient $B$
and the universal orbifold covering $N$ of $B$. The absence of conjugate points
implies flatness of $N$ in such a case. This flatness, on the other hand,
implies that the action is hyperpolar. 
Thus \tref{mainthm} generalizes Theorem A from \cite{LT}.

 As soon as one knows that the quotient $B=M/K$ is a Riemannian orbifold,
the statement of \tref{mainthm}  is a more or less direct consequence
 of the definition
of variational completeness. Thus the main part of the work  consists
in proving that variational completeness implies that the quotient $M/K$
is a Riemannian orbifold. In this part of the proof
\tref{theorem1} plays an essential role.

  As in \cite{LT}, our results generalize to singular 
Riemannian foliations.
Recall that a {\it transnormal system} $\mathcal F$ on a Riemannian
manifold $M$ is a decomposition of $M$ into smooth injectively immersed 
connected submanifolds, called {\it leaves},  such that geodesics 
emanating perpendicularly from one leaf stay perpendicularly to all
leaves. A transnormal system $\mathcal F$ is called a {\it singular Riemannian
foliation} if  there are smooth vector fields $X_i$
on $M$ such that for each point $p\in M$ the tangent space
$T_p  L (p)$ of the leaf $ L(p)$ through $p$ is given
as the span of the vectors $X_i (p) \in T_p M$. We refer to
\cite{Molino}, pp.185-216 and \cite{Wilk} for more on 
singular Riemannian foliations.
  Examples of singular Riemannian foliations are (regular) 
Riemannian foliations and  the orbit decomposition of an isometric group
action.  A singular Riemannian foliation $\mathcal F$ 
 will be called {\it closed}
if all of its leaves are closed in $M$, and it will be called
{\it locally closed at a point} $x\in M$  if for some neighborhood $U$ of $x$
the restriction of $\mathcal F$ to $U$ is closed (i.e., connected components
of the intersection of the leaves of $\mathcal F$ with $U$ are closed in $U$).
If $\mathcal F$ is locally closed at $x$ (if $\mathcal F$ is closed and $M$
is complete, respectively) then the local quotient 
(the global quotient, respectively)
$U/\mathcal F$ ($M/ \mathcal F$) is a 
well defined Alexandrov space of curvature
locally bounded from below  (\cite{BGP}, \cite{Bol},\cite{subm}).
Note, that closedness and completeness are only needed to assure that the 
quotient space is Hausdorff.

 For a singular Riemannian foliation $\mathcal F$ on $M$, denote by 
$M_0$ the set of all points  $z\in M$,
called {\it regular points}, whose leaves have maximal dimension.  
The set $M_0$ is open and dense in $M$, and the restriction of 
$\mathcal  F$ to $M_0$ is a (non-singular) Riemannian foliation.
  At each point $z\in M_0$, the local quotient of $M$ 
modulo   $\mathcal F$ is a smooth Riemannian manifold and we denote by 
$\bar \kappa (z)$ the maximum of all sectional curvatures at the image of $z$
in this quotient. The next result describes (non-) explosion 
of these curvatures as one approaches a boundary point $x$ of $M_0$.
The isotropy representation at 
$x$ used in  \tref{theorem1} has now to be replaced by the
infinitesimal singular Riemannian foliation $T_x \mathcal F$ on 
the tangent space $T_x M$ 
(see \cite{Molino}, pp.202-205 and Subsection \ref{infinit}).

 The following result generalizes \tref{theorem1}:

\begin{thm} \label{theorem2}
Let $\mathcal F$ be a singular Riemannian foliation on a Riemannian manifold
$M$. Let $x\in M$ be a point and let $T_x \mathcal F$ be the
infinitesimal singular Riemannian foliation induced by $\mathcal F$ on the tangent space $T_x M$.
Then the following are equivalent:
\begin{enumerate}
\item  $\limsup _{z\in M_0, z\to x}  \bar \kappa (z) < \infty$ ;

 \item  $\limsup _{z\in M_0, z\to x}  \bar \kappa (z) \cdot  d^2 (x,z) =0$;

 \item The singular Riemannian foliation $T_x \mathcal F$  is polar;

 \item $\mathcal F$ is locally closed at $x$ and a local quotient 
$U/\mathcal F$ of a neighborhood $U$ of $x$  is a Riemannian orbifold.
\end{enumerate}
\end{thm}

Let $\mathcal F$ be a singular Riemannian foliation on a Riemannian manifold
$M$. We will call $\mathcal F$ {\it infinitesimally polar } 
at the point $x\in M$
if the equivalent conditions of  \tref{theorem2} are fulfilled.
We call $\mathcal F$ infinitesimally polar if it is infinitesimally polar
at all points of $M$.  The class  of infinitesimally polar singular 
Riemannian foliation is a broad generalization of regular Riemannian 
foliations and singular Riemannian foliations with  sections. On the 
one hand, such singular Riemannian foliation  are well accessible by 
differential geometric methods, since their local quotients are 
finite(!) quotients of smooth Riemannian manifolds. 
On the other hand, each singular Riemannian foliation is infinitesimally 
polar on large parts of $M$ (compare \pref{smcodim}).
 We hope to discuss general features of infinitesimally polar 
singular Riemannian foliations  somewhere else.  Here we discuss 
a surprising characterization of infinitesimally polar singular
 Riemannian foliations
 and an application to general closed singular Riemannian foliations. 

 In order to state the characterization, recall that a geodesic on $M$
is called {\it horizontal} if it meets the leaves of $\mathcal F$ perpendicularly.
 We will
call such a geodesic $\gamma :[a,b] \to M$ {\it regular}, if $\gamma (a)$ and
$\gamma (b)$ are regular points of $\mathcal F$.  A regular horizontal 
geodesic 
intersects the singular strata of $\mathcal F$ only in finitely many points 
$a<t_1 <....<t_k <b$ (see \cref{reggeod}).
 We set $c(\gamma ) =c_1+....+c_k$, where $c_i$ is given by
$\dim (L(\gamma (a))) -\dim (L (\gamma (t_i)))$, and call this number
the {\it crossing number} of $\gamma$.   Then we have:

\begin{thm} \label{crosscont}
Let $M$ be a Riemannian manifold and let  $\mathcal F$ be a singular Riemannian
foliation on $M$. Let $G$ denote the space of all regular horizontal  geodesics
with the topology of pointwise convergence. Let $c:G \to \mathbb N$
be the crossing number function.  The function $c$ is continuous if and only
if $\mathcal F$ is infinitesimally polar.
\end{thm}

\begin{rem}
The crossing numbers $c(\gamma )$ were
 implicitly used in \cite{Bott} and \cite{BottS}
to calculate the indices of geodesics in compact symmetric spaces
and to deduce consequences about homology groups of their loop spaces. 
The equality between the index  and the crossing number holds true for 
{\it variationally complete} actions, or, more generally, for  singular
Riemannian foliations without horizontal conjugate points (see below).
Generalizing \cite{BottS}, one can deduce from this fact that the leaves
are taut submanifolds, and the crossing numbers can be used to compare
the topology of the leaves with the topology of the ambient 
space (cf. \cite{Nowak}).
 \end{rem}

Let now $\mathcal F$ be an arbitrary singular Riemannian foliation
on a complete Riemannian manifold $M$. Then  on the quotient $M/\mathcal F$
one can define a canonical ``quasi-geodesic flow''. Restricted to the
regular part $B_0$ of $M/\mathcal F$ this flow coincides with the geodesic flow. However,
at some time instances the flow may leave the regular part, and then
a local increase or decrease of volume could happen a priori. Using 
\tref{theorem2} we show that it does not happen. In case of a compact quotient,
the next result can be used to obtain almost recurrent horizontal geodesics.

\begin{thm} \label{quasi}
 Let $M$ be a complete Riemannian manifold and  let $\mathcal F$ be a closed
singular Riemannian foliation on $M$.  Then the projection of the horizontal
 geodesic flow on $M$ leaves the Liouville measure of $M/\mathcal F$
invariant.
\end{thm}

The notion of variational completeness generalizes to 
the setting of singular Riemannian foliation
as the notion of absence of horizontal conjugate points (cf. \cite{LT}).
Namely, let $\gamma :[a,b] \to M$ be a horizontal geodesic. An 
{\it $\mathcal F$-Jacobi field} along $\gamma$ is a variational field
through horizontal geodesics starting on the leaf $L(\gamma (a))$. An
{\it $\mathcal F$-vertical Jacobi field} along $\gamma$ is an 
$\mathcal F$-Jacobi field $J$
with $J(t) \in T_{\gamma (t)} L(\gamma (t))$ for all $t$. We say that 
{\it $\gamma$ has no horizontal conjugate points} if each $\mathcal F$-Jacobi
field $J$ with $J(t_0) \in T_{\gamma (t _0)} L(\gamma (t _0))$ for some 
$a<t_0< b$ is $\mathcal F$-vertical.  
We say that {\it $\mathcal F$ has no horizontal conjugate 
points} if no horizontal geodesics in $M$ have horizontal conjugate points.

The following result   generalizes    \tref{mainthm}:

\begin{thm} \label{maingen}
Let $M$ be a complete Riemannian manifold with a 
singular Riemannian foliation
$\mathcal F$.  If $\mathcal F$ has no horizontal conjugate points then it is 
infinitesimally polar. 
If $\mathcal F$ is closed then $\mathcal F$ has no horizontal conjugate points
if and only if  quotient $M/\mathcal F$ is a
good Riemannian orbifold without conjugate points.
\end{thm}

 The paper is structured as follows. In Section \ref{secprelim} we recall
some basic facts about singular Riemannian foliations and Riemannian orbifolds.
In Section \ref{secinfpol}, we prove \tref{theorem2}, 
its consequences \tref{theorem1}
and \pref{smcodim} that generalizes \cref{strata} to the case of 
singular Riemannian  foliations. In Section \ref{horflow} we discuss
the horizontal geodesic flow  and prove \tref{quasi}.  The main technical
observation of this section is the fact that the horizontal geodesic
flow in the total space defines a flow in a quotient, i.e., that two
projections of horizontal geodesics that coincide initially coincide
for their   life span.  This result was proved in \cite{subm}
and \cite{Bol} for the case of proper singular Riemannian foliations, and
in \cite{LT} and \cite{Nowak} for the case of singular Riemannian foliations
without conjugate points. An independent proof of this fact recently appeared
in \cite{ATequifoc}.  Finally, in Section \ref{secconj} we discuss various 
notions of conjugate points, prove stability of the absence of conjugate points
and  deduce \tref{maingen} and \tref{crosscont}.

We would like to thank Burkhard Wilking for several useful conversations.

\section{Preliminaries} \label{secprelim}
 Let $M$ be a Riemannian manifold and let $\mathcal F$ be a singular Riemannian
foliation on $M$.

\subsection{Distinguished tubular neighborhoods and 
infinitesimal foliations} \label{infinit}
Let $x\in M$ be a point. Then there is a small  open ball $P$  around $x$
in the leaf $L(x)$, a number $\epsilon >0$ and a neighborhood $O$
of $P$ in $M$, called a {\it distinguished tubular neighborhood at $x$}
such that the following holds true (\cite{Molino}, pp.192--193 and pp.202--205):

\begin{enumerate}
 
\item The foot point projection $F:O\to P$ is well defined;

\item $O$  is the image of the $\epsilon$-tube $N^{\epsilon} (P)$ in the 
normal bundle $N(P)$ of $P$ under the exponential map, and the map
$\exp :N^{\epsilon} (P) \to O$ is a diffeomorphism ;

\item For each real positive number $\lambda \leq 1$ the map 
$h_{\lambda} :O\to O$, given by  $h_{\lambda }  (\exp (v) )= \exp (\lambda v)$
for all $v\in N^{\epsilon } (P)$, preserves $\mathcal F$;

\item There is a diffeomorphism  $\phi $ of $O$ into the tangent space $T_x M$
  with $D_x \phi =\Id$ and a singular Riemannian foliation
 $T_x \mathcal F$ on  $T_x M$  that  coincides with
 $\phi _{\ast} \mathcal F  $ on $\phi (O)$ and such that 
$T_x \mathcal F$ is invariant under all rescalings 
$r_{\lambda} : T_x M\to T_x M, r_{\lambda } (v)=\lambda v$, for all 
$\lambda >0$. 
\end{enumerate}

\begin{rem}
In Section \ref{horflow} we will see that the restriction
of $\mathcal F$ to $O$ is 
invariant under the reflection $h_{-1}:O\to O $ at $P$.
\end{rem}

 The singular Riemannian foliation $T_x \mathcal F$ on the
tangent space $T_x M$ will be called the {\it infinitesimal 
singular Riemannian foliation
of $\mathcal F$ at the point $x$}.  
The infinitesimal foliation $T_x \mathcal F$ can be considered as a blow up
of $\mathcal F$ in the following sense.  Let $\phi$ be as above. Identify
$O$ with $\phi (O)$.
Set $O ^\lambda = r_{\lambda} (O)$ and define the Riemannian 
metric $g^{\lambda}$ on $O^{\lambda}$
as $g^{\lambda} =\lambda ^2 \cdot ( r_{\lambda}) _{\ast} g$. 
We have $\cup O^{\lambda} = T_x M$.
On compact subsets of $T_x M$ the 
blow up metrics $g^{\lambda}$ smoothly converge
to the flat metric $g_x$.  By construction, the restriction of
$T_x \mathcal F$ to $O^{\lambda}$ is a singular 
Riemannian foliation with respect to  $g^{\lambda}$.

\subsection{Local quotients}
We will continue to use the notations introduced above.
Note that $\mathcal F$ is locally closed at the point $x$ if and only
if the infinitesimal foliation $T_x \mathcal F$ is closed. In such
a case, the quotient $T_x M/T_x \mathcal F$ is a non-negatively curved
Alexandrov space and $\mathcal \phi (O) /T_x \mathcal F$
 is a ball around the origin (the leaf through $0$)
in this space.  The space $O/\mathcal F$  is an inner metric space
of curvature bounded below in the sense of Alexandrov. Moreover,
the space $T_x M/T_x \mathcal F$ is the tangent space to this Alexandrov
space at the leaf $L(x) \in  O/\mathcal F$.

Let us now assume that $M$ is complete and that $\mathcal F$ is closed.
Let $x\in M$ be given, and let $O$ be a small distinguished tubular
neighborhood of $x$, 
such that $O\cap L(x)=P$, in the notations of Subsection \ref{infinit}.
 Then $M/\mathcal F$
and $O/\mathcal F$ are spaces with curvature locally bounded below.
The embedding $O\to M$ induces an open map 
$i:O/\mathcal F \to  M/\mathcal F$. Since $\mathcal F$ is closed, the map
is finite-to-one and, by construction, the leaf $L(x) \in  M/\mathcal F$
has only one preimage in $O/\mathcal F$.  The map
$i$ preserves the lengths of all curves.  

 Let $U$ be an $\epsilon$-tube around the leaf $L$ with the same
$\epsilon$ as in the definition of $O$. Then $U$ is a union of
leaves of $\mathcal F$ and $U/\mathcal F$  is a neighborhood of
$L(x)$ in $M/\mathcal F$; see \cite{subm} or Section \ref{horflow}.
The global quotient $U/\mathcal F$ is mapped
by $i$ onto the local quotient $O/\mathcal F$.  To understand 
the map $i:O/\mathcal F \to  U/\mathcal F$, consider the universal covering
$\tilde U$ of $U$ with the lifted singular Riemannian foliation 
$\tilde  {\mathcal F}$ and the group of deck transformations $\Gamma$.
We have $\tilde U /\mathcal {\tilde F} =O/\mathcal F$. The action
of $\Gamma$ on $\tilde U$ induces an isometric action on 
$\tilde U/\mathcal {\tilde F}$ whose quotient is precisely $U/\mathcal F$.
Thus we deduce that a neighborhood of $L(x)$ in $M/\mathcal F$ is the 
quotient of the local quotient $O/\mathcal F$ by a finite group of isometries
of $O/\mathcal F$.

\subsection{Stratification} 
Let $\mathcal F$ again be a singular Riemannian foliation on the Riemannian
manifold $M$.
By the dimension, $ \dim  (\mathcal F)$, and the codimension of $\mathcal F$,
 $\codim  (\mathcal F, M)$,  we denote  the maximal 
dimension, respectively the minimal codimension of its leaves.  
For $s\leq \dim (\mathcal F)$ denote by $\Sigma _s$
the subset of all points $x\in M$  with $\dim (L(x))=s$. Then
$\Sigma _s$ is an embedded submanifold of $M$  and the restriction of
 $\mathcal F$ to $\Sigma _s$ is a Riemannian foliation 
(\cite{Molino}, pp.194-198). For a point
$x\in M$, we denote by $\Sigma ^x$ the connected component of $\Sigma _s$
through $x$, where $s =\dim (L(x))$.   We call the decomposition of $M$
into the manifolds $\Sigma ^x$ the {\it canonical stratification} of $M$.

 The subset $\Sigma _{\dim (\mathcal F)}$ is open, dense and connected
in $M$. It is the regular stratum $M_0$ of $M$. All other singular strata 
$\Sigma ^x$ have codimension at least $2$ in $M$.    The {\it quotient 
codimension of a stratum} $\Sigma  ^x$ is defined to
be $\codim (\mathcal F ,M) - \codim (\mathcal F , \Sigma ^x)$.

  If the singular Riemannian foliation $\mathcal F$ is locally
closed at $x$, then a local quotient $O/\mathcal F$ is a space
stratified by smooth Riemannian orbifolds 
 $\Sigma ^z / \mathcal F$  and the quotient
codimension of the stratum $\Sigma ^x$ is just the codimension
of the quotient $\Sigma  ^x /\mathcal F$ in the whole quotient 
$O/\mathcal F$.

  Let a point $x\in M$ be fixed. Then the tangent space $T_x M$
decomposes as $T_x M =T_x \Sigma ^x \oplus N_x \Sigma ^x$,
where $N_x \Sigma ^x$ is the normal space in $M$ to $\Sigma _x$.
The infinitesimal Riemannian foliation $T_x \mathcal F$ on $T_x M$
is the direct  product of the foliation of $T_x \Sigma ^x$ by affine subspaces
parallel to $T_x L(x)$ and a singular Riemannian foliation 
$\widetilde   {T_x \mathcal F}$ on $N_x \Sigma ^x$. 
This last (the main) part $\widetilde   {T_x \mathcal F}$ is invariant 
under positive homotheties of $N_x \Sigma ^x$ and the only $0$-dimensional leaf
of $\widetilde   {T_x \mathcal F}$ is the origin $\{ 0\}$.

 The quotient codimension of the stratum $\Sigma ^x$ is the codimension
of the singular Riemannian foliation  $\widetilde   {T_x \mathcal F}$ on
the Euclidean space  $N_x \Sigma ^x$.

\subsection{Riemannian orbifolds} \label{rimorb}
We refer to \cite{Bridson}, pp.584-619,
 for a more advanced and refined study of Riemannian
orbifolds. 
A metric space $X$ is called a {\it good
Riemannian  orbifold } if $X$ is isometric to $M/\Gamma$, where $M$
is a smooth Riemannian manifold and $\Gamma $ a discrete group of isometries.

 A point $x$ in a metric space $X$ is called an {\it orbifold point} if $x$
has a neighborhood $U$ that is a good Riemannian orbifold.  The set $O$
of all orbifold points in $X$ is open. We call $X$ a {\it Riemannian 
orbifold } 
if $X=O$ holds. Note that the quotient of a Riemannian orbifold by a 
finite group of isometries is again a Riemannian orbifold.

   Let $B$ be a Riemannian orbifold. Then locally $B$ is a finite
isometric quotient of a smooth Riemannian manifold. Since geodesics,
tangent spaces and the Liouville measure on the unit tangent bundle
of Riemannian manifolds are invariant under isometries, one
gets corresponding notions on $B$.  Namely, the ``unit tangent bundle''
$UB$   being a disjoint union  $UB=\cup _{y\in B} S_y B$ of spaces
of directions $S_y B$, locally being a finite quotient of the unit tangent 
bundle of a ``covering'' Riemannian manifold.
This unit tangent bundle comes along with the foot point projection
$p:UB\to B$, a locally compact (quotient) topology, a local geodesic
flow $\phi$ and the Liouville  measure $\mu$, that is a Borel measure on $UB$.
The local flow $\phi _t$ preserves the Liouville measure, whenever
it is defined, since this is the case for Riemannian manifolds and since
the Liouville measure and the local geodesic flow are preserved under 
isometries.  
For $v\in UB$,
we set $\eta _v (t) = p (\phi _t (v))$, i.e. the curve $\eta _v$
is locally the image of a geodesic in a Riemannian manifold under
the quotient map. We call $\eta_v$ the {\it orbifold-geodesic} in the direction
of $v$.

   A Riemannian orbifold $B$ is stratified by 
Riemannian manifolds with a unique maximal
stratum $B_0$ that is open and dense in $B$.  The unit tangent
bundle  $UB_0$ is an open and dense subset of $UB$, that has full measure
with respect to $\mu$. Moreover, the set $U'$ of all vectors $v\in UB_0$
such that the orbifold-geodesic $\eta _v$ does not cross strata of 
codimension $\geq 2$ is of full measure in $UB$.

For each orbifold-geodesic $\eta _v$, the curvature endomorphism  along
$\eta _v$ is well defined. Therefore,  the notions of Jacobi fields
and conjugate points are also well-defined. Let us now assume that $B$ is complete as a metric space.
Then each orbifold-geodesic is defined on $\mathbb R$ and the local 
geodesic flow is a global flow.
Take  a regular point $x\in B_0 \subset B$. Consider the orbifold exponential
map $\exp :T_x B \to B$ given by $\exp (tv) =\eta _v (t)$, for a unit vector
$v\in T_x$.  This map (since defined in metric terms) factors over
local branched covers of $B$, i.e. for each $w\in T_x B$ there is a
finite quotient $N/\Gamma _w =O \subset B$, with $\exp (w)\in O$,
such that $\exp$  lifts on a neighborhood of $w$ to a smooth map to $N$.
The vector  $w=tv$ is a conjugate vector along the geodesic $\eta_v$,
if and only if this lift has a non-injective differential at $w$.

If no vector in $T_x B$ is a conjugate vector of $x$, i.e., if $x$ has no
conjugate points, then one can pull back the metric from $B$
(in fact from the local covers $N$) to a Riemannian metric $\tilde g$ on 
$T_x B$.  The space $T_x B$ with this Riemannian metric $\tilde g$
is a complete Riemannian manifold and the map $\exp :T_x B\to B$
becomes an arclength preserving orbifold covering 
(as in the Theorem of Cartan-Hadamard).
Thus $T_x B$ is the universal orbifold covering of $B$
and  we get $B=(T_x B,\tilde g) /\Gamma$ for some group $\Gamma$ of isometries
of $T_x B$. Hence $B$ is a good orbifold in this case. Thus we have 
the following observation that is implicitly contained in the proof 
of the developability results stated in \cite{Bridson}, p.603.

\begin{lem} \label{developp}
 If $B$ is a complete Riemannian orbifold without conjugate points
then $B$ is a good orbifold. More precisely, there is a complete 
Riemannian manifold $N$ without conjugate points and a discrete group $\Gamma$
of isometries of $N$ such that $B=N/\Gamma$.
\end{lem}

\section{Infinitesimal polarity}  \label{secinfpol}

\subsection{Horizontal sections}
Let $\mathcal F$ be a singular Riemannian foliation on a Riemannian manifold
$M$. A {\it global (local) horizontal section through} $x$ is a 
smooth immersed 
submanifold $ N$ in  $ M$ through $x$ 
that intersects all leaves of $\mathcal F$ 
(all leaves in a neighborhood of $x$), such that all intersections are
orthogonal. We say that $\mathcal F$  is  {\it polar} if there are 
global horizontal
sections through every point $x\in M$. It is called {\it hyperpolar} 
if all these
sections are flat.  Recall that all local sections are totally geodesic;
thus each polar foliation of $\mathbb R^n$ is hyperpolar. Finally,
each hyperpolar foliation of $\mathbb R ^n$ is closed.
We refer to \cite{Boualem},\cite{Alexandrino},\cite{Alexandrino1}
 for more on singular Riemannian foliation with sections.

\begin{rem}
In many important cases, like in  Euclidean  or  symmetric spaces, 
polar  singular Riemannian foliation  have been objects of an extensive
study, where  they are better known  as isoparametric
foliations. See, for instance, \cite{Palais}.
\end{rem}

\subsection{Non-explosion of curvature} Now we can start with the

{\it Proof of  \tref{theorem2}.}
 The implications
$(4) \Longrightarrow (1) \Longrightarrow (2)$ are clear. 

 Assume $(2)$. We are going to use the notations introduced in Subsection
\ref{infinit}.  Let $z\in T_x M$ be a $T_x \mathcal F$-regular point.
For all $t$ with  $z\in O^t$ we denote by $\bar \kappa ^t (z)$
the supremum of all sectional curvatures at the local projection of $z$ in 
$(O^t,g^t)/T_x \mathcal F$. Since $g^t$ smoothly converge to $g_x$
for $t\to \infty$ the horizontal curvatures satisfy $\lim _{t\to \infty} 
(\bar \kappa ^t(z) ) =\bar \kappa _x (z)$. On the other hand, the assumption
$(2)$ implies $\lim  _{t\to \infty} (\bar \kappa ^t(z) ) =0$. Thus we deduce
that the local quotient of $(T_x M,g_x)$ modulo $T_x \mathcal F$  is flat
at all regular points. Due to the flatness of $(T_x M,g_x)$, this implies the 
vanishing of the O'Neill tensor of the Riemannian foliations  $T_x \mathcal F$
on the regular part. Thus the horizontal distribution on the regular part
of $T_x \mathcal F$ is integrable and we deduce  that
$T_x \mathcal F$ is hyperpolar (\cite{Alexandrino1}).

The main implication $(3)\Longrightarrow (4)$ is  more subtle. Let 
$T_x \mathcal F$ be hyperpolar and let $N$ be a 
horizontal section through $0$. There is  a finite group $\Gamma$ 
of isometries of $N$, called the{\it  Weyl group} of $N$ (cf. \cite{Palais}), 
with $\Gamma (0)=0$ such that $N/\Gamma$ is isometric
to $T_x M /T_x \mathcal F$. Let $O$ be again as in Subsection \ref{infinit}
and let us again identify it with $\phi (O)$.  Let $N_0$ be a small
ball in $N$ around $0$ that is  contained in $O$. 

 In general, one cannot expect that $N$ is a $g$-horizontal section of 
$\mathcal F$ nor must $\mathcal F$ have any horizontal local sections.
The idea is  to define a new ``horizontal'' metric on $N_0$,
that is invariant under $\Gamma$ and such that $N_0 /\Gamma$ becomes
isometric to a neighborhood of $x$ in $O/\mathcal F$.

 For a point $z\in N_0$, we denote by $\tilde V (z)$ the orthogonal complement
of $T_z N$ with respect to the (constant) metric $g_x$.  
Let $\tilde H(z)$ denote
the orthogonal complement of $\tilde V(z)$ with respect to the original metric
$g$. 
Then $\tilde H (z)$ depends smoothly on $z$.
 Moreover, each space $\tilde H(z)$ is contained in the $g$-orthogonal
complement of $T_z L (z)$. By dimensional reasons, $\tilde H(z)$ is the
orthogonal complement of $T_z L(z)$ at all regular points.

   Define  a  Riemannian metric $\tilde g$ on $N_0$  by 
$\tilde g _z (v,w)=  g_z (P^z(v), P^z(w))$, where $P^z$ is the 
orthogonal projection
to $\tilde H (z)$ with respect to $g_z$. 
Denote by $\tilde N$ the manifold $N_0$ with the inner metric 
defined by the Riemannian metric $\tilde g$. The projection
$p:\tilde N \to O /\mathcal F$ preserves lengths of all curves contained in 
the set of regular points of $N$. On the other hand, $p$ is invariant
under the group of diffeomorphisms $\Gamma$ of $\tilde N$.
Thus each $k\in  \Gamma$ preserves lengths of all
curves contained in the regular part of $N$.
 By continuity,  each $k\in  \Gamma $
is an isometry of $\tilde N$. 
Moreover, the induced map $\tilde N /\Gamma  \to O/\mathcal F$
is an isometric embedding. This proves $(4)$ and finishes the 
proof of \tref{theorem2}. 
\qed

  Let us now assume that $M$ is complete and that $\mathcal F$ is closed.
 Let $x\in M$ be given and let $O$ be a small distinguished 
tubular neighborhood around $x$.  If $\mathcal F$  is infinitesimally polar,
then $O/ \mathcal F$ is a Riemannian orbifold. Therefore,
the image of $O$ in $M/\mathcal F$, that is a finite quotient
of $O/ \mathcal F$, is a Riemannian orbifold. On the other hand, if
$L(x)\in  M/\mathcal F$ is an orbifold point then the regular
part of $T_x M/ T_x \mathcal F$ is flat and we get that $T_x \mathcal F$
is hyperpolar. Thus $L(x)$ is an orbifold point
of the global quotient $M/\mathcal F$ if and only if $\mathcal F$ is 
infinitesimally polar at $x$. Thus  \tref{theorem2} implies
\tref{theorem1}.

\subsection{Small codimensions}
Let $\mathcal F$ be a singular Riemannian foliation
on the Euclidean space $\mathbb R^n$, that is 
invariant under positive rescalings and satisfies  $L(0)= \{ 0\}$.
Then all leaves of $\mathcal F$
are contained in concentric spheres around $0$, and $\mathcal F$ is the cone
over the restriction of $\mathcal F$ to the unit sphere $\mathbb S^{n-1}$.
Note, that $\mathcal F$ is polar on $\mathbb R^n$ if and only if its 
restriction to $\mathbb S^{n-1}$ is polar.
We have $\codim (\mathcal F, \mathbb R^n) =
\codim (\mathcal F , \mathbb S^{n-1} ) +1$. 
Finally, each singular Riemannian foliation of codimension $1$
in a complete Riemannian manifold is polar. Thus each scaling invariant
singular Riemannian foliation on  
$\mathbb R^n$ of codimension $\leq 2$ is polar.
This proves the following result generalizing \cref{strata}:

\begin{prop} \label{smcodim}
 Let $\mathcal F$ be a singular Riemannian foliation on a Riemannian 
manifold $M$.  Let $x\in M$ be a point with stratum $\Sigma ^x$
of quotient codimension $\leq 2$. Then $\mathcal F$ is infinitesimally
polar at $x$.
\end{prop}

\section{Horizontal exponential map} \label{horflow}
\subsection{Horizontal vectors}
Consider the subset $D$ of the unit tangent bundle $UM$ that consists
of all starting vectors $v$ of horizontal geodesics $\gamma _v$. 
The set $D$ is closed and invariant under the local geodesic flow $\phi _t$,
whenever it is defined.  By $p:D\to M$ we denote the foot point projection.
 
We now discuss  a preliminary  stratification of $D$; later we will derive
a more natural stratification adapted to the geodesic flow.
 For each stratum $\Sigma ^x$ of $M$, the preimage $D^x = p^{-1} (\Sigma ^x) 
\subset D$ is a smooth submanifold of $UM$ (and in fact  a sphere bundle 
over $\Sigma ^x$) of dimension $\dim (M) -1 + \dim (\Sigma ^x) -\dim (L(x))$.
Thus $D$ is stratified by smooth submanifolds $D^x$ of $UM$. The main stratum
$D_0$ is open and dense in  $D$ and the codimension of the stratum $D^x$
(i.e., $\dim (D_0) - \dim (D^x)$) coincides with the quotient codimension
of the stratum $\Sigma ^x$ of $M$.

 Let $M_i$ be the open subset of all points $x\in M$ such that the quotient
codimension of $\Sigma ^x$ is at most $i$. Set $D_i = p^{-1} (M_i) \subset D$.
Then $D_i$ is the union of all strata in the above stratification of $D$ that
have codimension at most $i$.  Since the geodesic flow $\phi _t$ is locally
Lipschitz, the shadow of $D\setminus D_i$  under $\phi$  (i.e.,
the set of all directions $v\in D$, such that $\gamma _v$ intersects 
$M \setminus M_i$)  has Hausdorff dimension at most $\dim (D) -i$.  
Moreover, for each fixed $t$, the 
Hausdorff dimension of $\phi _t (D\setminus D_i)$ is at most $\dim (D) -i-1$. 
In particular, we deduce:

\begin{lem} \label{outside}
There is a subset $D'$ of full measure in the manifold $D_0$, such that
for all $v\in D'$ the whole geodesic $\gamma _v$ is contained in $M_1$.
\end{lem}

\subsection{Horizontal geodesics in the nice part}   \label{nicepart}
Let $x$ be a regular point ($x\in M_0$). Then 
a local quotient of $O /\mathcal F$ around $x$ 
is a smooth Riemannian manifold, and horizontal geodesics in $O$
are projected to geodesics in $O/\mathcal F$. In particular, two
such projections coincide, if they coincide initially.

Let now $x$ be a point, such that $\Sigma ^x$ has  quotient codimension $1$. 
Then the infinitesimal
quotient $T_x M /T_x \mathcal F$ is isometric to $\mathbb R ^{q-1} \times
[0, \infty )$ and (due to \tref{theorem2})
the local quotient  $O/\mathcal F$ is a Riemannian
orbifold of the form $N/ \mathbb Z_2$, where $\mathbb Z_2$ acts
as an isometric  reflection at a 
totally geodesic hypersurface of a smooth Riemannian
manifold $N$.  A horizontal geodesic $\gamma$ in $O$  
is either completely
contained in the regular part of $O$, or it is completely contained
in the singular stratum $\Sigma ^x$, or it intersects $\Sigma ^x$
in precisely one point. 
In the first two cases the image of $\gamma$ is a geodesic in $O /\mathcal F$.
In the last case the image  is the concatenation of two 
geodesics that meet at the
boundary of $O /\mathcal F = N/\mathbb Z_2$ and satisfy the reflection
law. In any case,
such a projection is an orbifold-geodesic in $N /\mathbb Z_2$.  
In particular, two  projections of horizontal 
geodesics coincide if they coincide initially.     Thus we have shown:

\begin{lem} \label{lnicepart}
If $M$ does not have any strata of quotient codimension $\geq 2$, then
in each local quotient $O/\mathcal F$  projections of horizontal 
geodesics coincide if they coincide initially.
\end{lem}

\subsection{Equivalence relation} 
We are going to define   a natural equivalence relation 
$\mathcal R$ on $D$ that identifies 
two directions $v,w \in D$ if the corresponding geodesics have equal images
in $M/\mathcal F$.  

 To be more precise, let first $\mathcal  L$ be a leaf of $\mathcal F$.
There is a small (not necessarily tubular) neighborhood $U$ of 
the zero section of the normal bundle
$N(\mathcal L)$ of $\mathcal L$ with the following properties (see Subsection
\ref{infinit}).
The exponential map  restricted to $U$ is a local diffeomorphism.
The set $U$ is pointwise star-shaped, i.e. it is invariant under the
maps $h_{\lambda } (v) =\lambda v$,  $0< \lambda \leq 1$. Finally,
the pull back $\exp ^{\ast} (\mathcal F)$ is invariant under all 
$h_{\lambda}$. Thus there is a unique singular foliation
on $N(\mathcal L)$ invariant under all $h_{\lambda}$, $0<\lambda < \infty$,
 that coincides with 
$\exp ^{\ast} (\mathcal F )$ on $U$. 

We will call two vectors $v,w \in D$ equivalent if $v$ and $w$
are normal vectors  to the same leaf $\mathcal L$ of $\mathcal F$,
and if they are in the same leaf of the singular foliation on the normal 
bundle $N(\mathcal L)$ described above.    Equivalently, $v$ and $w$ are in
the same equivalence class if and only if there is a smooth 
(or piecewise smooth) curve $\eta $ connecting $v$ and $w$ in $D$
and a small positive number $\epsilon$, such that
the leaf through $\gamma _{\eta (s)} (t)$ does not depend on $s$, 
for all $0\leq t < \epsilon$. The last condition just means, that
for all $0\leq t < \epsilon$ the curve $\eta _t  (s)= \phi _t (\eta (s))$
is contained in the normal bundle to some leaf $\mathcal L _t$
 of $\mathcal F$. 

We will denote the equivalence relation by $\mathcal R$. 
By $\mathcal R(v)$ we will denote the equivalence
class of $v$.  Note that the restriction 
of $\mathcal R$ to the manifold $D_0$ is given by leaves of a smooth foliation.

We are going to prove the invariance of $\mathcal R$ under the geodesic flow
(cf. \cite{ATequifoc} for an alternative proof and \cite{subm}, \cite{Bol}
and \cite{Nowak} for some special cases).
\begin{prop} \label{invgeod}
Let $\eta :[a,b] \to D$  be a curve  in an equivalence class $\mathcal R (v)$.
If $\phi _t (\eta )$ is defined for some $t >0$, then $\phi _t (\eta )$
is contained in an equivalence class of $\mathcal R$. Moreover, 
$\mathcal R$ is invariant under the reversion $-\Id :D\to D$, given by
$-\Id (v) =-v$.
\end{prop}

{\it Proof.}
The equivalence classes of $\mathcal R$ are smooth injectively immersed 
submanifolds of the unit tangent bundle of $M$. Denote by 
$\tilde {\mathcal R } (v)$ the tangent space to the vector $v\in D$ of its 
equivalence class $\mathcal R (v)$. The claim can now be restated as follows:
The local flow $\mathcal \phi _t$ and the reversion $- \Id$  leave the
``distribution'' $\tilde {\mathcal R}$ invariant, i.e. for all
$v\in D$ we have $\tilde {\mathcal R} (- v)= 
(- \Id )_{\ast}  (\tilde {\mathcal R} (v))$ and $(\phi _t )_{\ast} 
(\tilde {\mathcal R} (v)) =\tilde {\mathcal R} (\phi _t (v))$, for all
$t$ such that $\phi _t (v)$ is defined.

  Given an open subset $V$ of $M$, the restriction of $D$, 
$\tilde {\mathcal R}$ and the flow $\mathcal \phi$ to the unit tangent
bundle of $V$ coincides with the corresponding  objects for the 
restriction of $\mathcal F$ to $V$. Thus our claim is local on $M$. 

 Due to  \lref{lnicepart}, the claim is true if in $M$ there are
no strata of quotient codimension $\geq 2$. Thus for all $v$, such that
the geodesic $\gamma _v :[0,t] \to M$ is contained in $M_1$, we have
 $(\phi _t )_{\ast} 
(\tilde {\mathcal R} (v)) =\tilde {\mathcal R} (\phi _t (v))$.

 In particular, this is true for all $v\in D' \subset D_0$ from 
\lref{outside} and all $t$, such that $\phi _t (v)$ is defined.
Notice that $\tilde {\mathcal R}$ is a smooth foliation on $D_0$, $\phi $
is smooth and $D'$ is dense in $D_0$. Therefore,  for all $v\in D_0$
and all $t$ with $\phi _t (v) \in D_0$ we must have
$(\phi _t )_{\ast} 
(\tilde {\mathcal R} (v)) =\tilde {\mathcal R} (\phi _t (v))$.  

\begin{rem}
Continuity arguments could be used to   finish the proof at this point if
all leaves were assumed to be  closed. 
\end{rem}

 Let $x\in M$ be a point and let the plaque  $P\subset L(x)$, the number
 $\epsilon >0$ and
 a distinguished tubular neighborhood $O$
at $x$ be chosen as in Subsection \ref{infinit}.   For a unit normal
vector $v$ to the plaque $P$, we get from the definition of $O$ and 
$\mathcal R$, that $(\phi _t )_{\ast} 
(\tilde {\mathcal R} (v)) =\tilde {\mathcal R} (\phi _t (v))$ and
$(\phi _t )_{\ast} 
(\tilde {\mathcal R} (-v)) =\tilde {\mathcal R} (\phi _t (-v))$, for
all $0\leq t <\epsilon$.  Moreover, 
$(\phi _t )_{\ast} 
(\tilde {\mathcal R} (v)) =\tilde {\mathcal R} (\phi _t (v))$ for all
$-\epsilon < t <0$  if and only if $(- \Id)_{\ast} (\tilde {\mathcal R} (v))=
\tilde {\mathcal R} (-v)$.

As we have seen, 
 $\phi _t$ leaves $\tilde {\mathcal R}$ on the regular part $D_0$ invariant,
therefore $(- \Id)_{\ast} (\tilde {\mathcal R} (v))=
\tilde {\mathcal R} (-v)$, for all normal $v$ to $P$, such that $\exp (tv)$
and $\exp (-tv)$ are in $M_0$, for some (and hence all) $0<t<\epsilon$.

 Consider the diffeomorphism $I:O\to O$ (reflection at $P$; in terms
of Subsection \ref{infinit} it is just $h_{-1}$), defined
by $I(\exp (tv))=\exp (-tv)$ for unit normal vectors $v$ to $P$ and
$0\leq t<\epsilon$.  By definition, $(-\Id )_{\ast} (\tilde {\mathcal R} (v))=
\tilde {\mathcal R} (-v)$ for a unit normal vector $v$ to $P$ if and 
only if $I$ preserves $\mathcal F$ at $\exp (\frac \epsilon 2 v)$.
 By the observation above,
$I$ preserves $\mathcal F$ on the open dense subset $M_0 \cap O \cap 
I(M_0 \cap O)$ of $O$. But a singular Riemannian foliation is uniquely defined 
by its restriction to an open dense subset; see \lref{uniquenessfol} below.
 We deduce $I_{\ast} (\mathcal F)
=\mathcal F$. Thus we have shown the invariance of $\tilde {\mathcal R}$
under the reversion $- \Id$.

 For each vector $v\in D$, we can now take its foot point $x$ and a 
distinguished tubular neighborhood of $x$ and deduce that
$(\phi _t )_{\ast} 
(\tilde {\mathcal R} (v)) =\tilde {\mathcal R} (\phi _t (v))$, for
all $-\epsilon < t < \epsilon $, where $\epsilon =\epsilon (x)$ is chosen as
in Subsection \ref{infinit}. Covering an arbitrary geodesic 
$\gamma _v :[0,t] \to M$ by finitely many distinguished tubular neighborhoods,
we deduce $(\phi _t )_{\ast} 
(\tilde {\mathcal R} (v)) =\tilde {\mathcal R} (\phi _t (v))$.

 This finishes the proof of \pref{invgeod}.
\qed

 In the proof above  we used the following:
\begin{lem} \label{uniquenessfol}
Let $M$ be a manifold. For $i=1,2$, let $g_i$ be a Riemannian metric on $M$
and let $\mathcal F_i$ be a singular Riemannian foliation on $M$ with respect
to $g_i$. If $\mathcal F_i$ coincide on an open and  dense subset $U$ of $M$ 
then they coincide on all of $M$.
\end{lem}

{\it Proof.}
 Choose an arbitrary point $x\in M$. The claim is local, thus restricting to a
small relatively compact neighborhood $O$ of $x$, 
we may assume that the leaves  $L_1 (x)$ and $L_2(x)$ are closed.  Then for 
each sequence $x_n \to x$ the leaves $L_i (x_n)$ (or, equivalently, their 
closures) converge in the Gromov-Hausdorff topology to the  leaf $L_i (x)$.

 Thus it is enough to prove that $L_1 (x)= L_2 (x)$ if $x$ is a regular
point of $\mathcal F_1$. Hence we may assume that $\mathcal F_1$
is a regular foliation. By continuity the leaves of $\mathcal F_2$ are
contained in the leaves of $\mathcal F_ 1$  
in such a case. Thus for each 
$\mathcal F_2$-regular point $y$ we get $L_1 (y)= L_2 (y)$. Then the above
limiting argument shows that leaves of $\mathcal F_1$ and $\mathcal F_2$
through all points coincide.
\qed

\subsection{Natural stratification of the space of 
horizontal geodesics} \label{natdecomp}

Let $v$ be a horizontal vector and let $\gamma =\gamma _v$ be the
horizontal geodesic with $\gamma _v ' (0) =v$. For
$t$ in the interval of definition of $\gamma $, we let $l(t)$ be
the dimension of the leaf $L(\gamma (t))$.   Due to the semi-continuity 
of the leaf dimension, we have $\liminf _{t_i \to t} l(t_i) \geq l(t)$.

 Let $t$ be fixed and consider a small distinguished neighborhood $O$ of
$\gamma (t)$ as in Subsection \ref{infinit}.  Then $h_s(\gamma (t+\rho ))=
\gamma (t+s\rho )$, for $-1\leq s \leq 1$.  For $0<s\leq 1$,
$h_s$ preserves $\mathcal F$.   On the other hand, in the course of
the proof of \pref{invgeod}, we have seen that $h_{-1}$ preserves $\mathcal F$
as well.  Thus $h_s:O\to O$ preserves $\mathcal F$ 
for all $s\in  [-1,1] \setminus \{ 0\}$. In  particular,
$l(t+\rho )$ does not depend on $\rho$ for 
$\rho \in     [-\epsilon , \epsilon ] \setminus \{ 0\}$.

 This shows that $l(t)$  is equal to the 
constant $d(\gamma ) := \max (l(t))$ for
all but discretely many $t$. We summarize our observations:

\begin{lem} \label{geoddim}
Let $\gamma$ be a horizontal geodesic in $M$. Let $d(\gamma )$ denote
the maximal dimension of the leaves $L(\gamma (t))$. 
Then for all but discretely many $t$, the leaf $L(\gamma (t))$
has dimension $d(\gamma )$.   
\end{lem}

 In the case of maximal dimension we get:

\begin{cor} \label{reggeod}
A compact horizontal geodesic that contains a 
regular point is contained in  the
set of regular points, with exception of at most 
finitely many points. 
\end{cor}

\begin{rem}
In fact, we have shown a slightly more general statement than \lref{geoddim}.
Namely, for each horizontal geodesic $\gamma$ and for all but discretely many
times $t$, the infinitesimal foliation $T_{\gamma (t)} \mathcal F$
does not depend on $t$.
\end{rem}

 By definition, the function  $d:D\to \mathbb N$, given by 
$d(v) =d(\gamma _v )$, is invariant under the local geodesic flow
$\mathcal \phi$  under the multiplication by $-1$ and 
under the equivalence relation $\mathcal R$.
From the above we deduce, that $d(v)$ is the dimension of the leaf
$L(\gamma _v (\epsilon ))$ for small  positive $\epsilon$.

 Set $D^i :=d^{-1} (i) \subset D$.  Then $D$ is decomposed into
the disjoint union of the sets $D^i$. Due to the  semi-continuity of leaf
dimensions, the closure of $D^i$ is contained in the union of $D^j$, 
$j\leq i$. 
We claim that $D^i$ is a  smooth  submanifold 
of the unit tangent bundle $UM$ of $M$. To see this, let $v\in D^i$
be given. Then $\phi _{\epsilon} (v)$ is a horizontal vector
of the restricted Riemannian foliation $\mathcal F$ on the submanifold
$\Sigma ^{\gamma _v (\epsilon)}$ of $M$. The space $\mathcal H$ of unit
horizontal vectors of the restriction of $\mathcal F$ to
the manifold $\Sigma ^{\gamma _v (\epsilon)}$ is a smooth
submanifold of the unit tangent bundle of  $\Sigma ^{\gamma _v (\epsilon)}$.
By definition, the diffeomorphism $\phi _{-\epsilon}$ sends a neighborhood of 
$\phi_{\epsilon } (v)$ in $\mathcal H$ to a neighborhood of $v$ in $D^i$.
This shows that $D^i$ are submanifolds of $UM$.

  If $\mathcal F$ is a Riemannian foliation  then the function
$d:D\to \mathbb N$ defined above is constant. For each vector
$v\in D$, we have $d(v)=\dim (L(p (v)))=\dim (\mathcal R(v))$.
Moreover, the equivalence classes of $\mathcal R$ are 
leaves of a foliation on the manifold $D$ in this case.
Finally, the equivalence classes of $\mathcal R$ are closed
if $\mathcal F$ is a closed Riemannian foliation.

 For a general singular Riemannian foliation, the observation above 
shows that each small open subset of $D^i$, for each $i$, can be moved
by the geodesic flow to an open part of the set of horizontal
vectors of a smooth Riemannian foliation (restriction of $\mathcal F$
to a stratum).  Thus we arrive at the following:

\begin{prop}
 For each $i\in D$, the subset $D^i$ of all vectors $v\in V$
with $\dim (\mathcal R (v)) =i$ is a submanifold of the unit tangent 
bundle, which  is invariant under the local geodesic flow $\phi$. 
The equivalence
classes of $\mathcal R$ on $D^i$ are leaves of a smooth foliation. This 
foliation has closed leaves, if the leaves of $\mathcal F$
are closed. For each $v\in D^i$, we have $d(\gamma _v) =i$.
\end{prop}

\subsection{Vertical Jacobi fields} \label{vertjac}
Let $v$ be a horizontal vector and $\gamma =\gamma _v$ the geodesic
in the direction $v$. Consider the leaf $\mathcal R (v)$ and a small
neighborhood  $V$ of $v$ in $\mathcal R (v)$.  For each $\bar v$
in $V$ consider the horizontal geodesic $\gamma _{\bar v}$. Due
to \pref{invgeod}, for each compact interval $I$ of definition of $\gamma$,
we may choose $V$ so small  that for all $t\in I$
and all $\bar v\in V$ we have  $\gamma _{\bar v} (t)  \in L (\gamma (t))$.

 The space $W^{\gamma }$ of variational fields  through
geodesics $\gamma _{\bar v}$ is 
a vector space of Jacobi fields along $\gamma $ of dimension
 $\dim (\mathcal R(v)) =d(v)$.  Due to \pref{invgeod},
$W^{\gamma }$ does not depend on the starting point of $\gamma$. Moreover,
we have $W^{\gamma} (t) := \{ J(t) | J\in W^{\gamma} \} = 
T_{\gamma (t)} L(\gamma (t))$, for all $t$ in the interval of definition of 
$\gamma$.

\subsection{Invariance of the Liouville measure}
Let now $M$ be complete and let $\mathcal F$ be closed.  
Consider the space of horizontal vectors $D$
and its decomposition $D=\cup D^i$ discussed in Subsection \ref{natdecomp}.
We have $v\in D^i$ if and only if $\dim L(\gamma _v (t)) =i$,
for all but discretely many times $t$.   The relation
$v\in \mathcal R(w)$ is equivalent to $L(\gamma _v (t))=L( \gamma _w (t))$,
for  all  $t \in \mathbb R$.

 Denote by $G$ the space of equivalence classes $D/\mathcal R$.  
 We consider
$G$ with the induced quotient topology 
(that is Hausdorff and locally compact in our case).
The decomposition of $D$ induces 
a decomposition of $G$ as $G=\cup G ^i$.  The flow
$\mathcal \phi$ descends to a flow on $G$.  Denote by $D^+$ and $G^+$
 the maximal stratum of $D$ and its projection to $G$, i.e., the set
of all  horizontal geodesics that contain at least one regular point.

 The subspace $D_0 \subset D^+$ of all horizontal vectors with
regular starting point is of full measure in the manifold $D^+$
(since the complement of $D_0$ is a countable union of submanifolds
of positive codimension).  The subspace $D' \subset D_0$ defined 
in \lref{outside} is of full measure in $D^+$, invariant under 
$\phi$ and saturated under $\mathcal R$.    

 For each $v\in D'$, the singular Riemannian foliation $\mathcal F$
is infinitesimally polar at all points $\gamma _v (t)$. Moreover,
due to Subsection \ref{nicepart}, $\gamma _v$ projects to an orbifold
geodesic in each local quotient $O/\mathcal F$.  Thus the image of
$\gamma _v$ in $M/\mathcal F$ is contained in the set $B$
of orbifold points of $M/\mathcal F$ and this image is an orbifold
geodesic in the orbifold $B$.

 We set $G' =D' /\mathcal R$. Identify $G'$ with a subset of the unit
tangent bundle $UB$ of the orbifold $B$.  The argument above shows that
on $G'$ the flow $\phi$ coincides with the  orbifold-geodesic 
flow of the orbifold $B$.

 Define the measure $\mu$ on $G$ by setting $\mu (G\setminus UB)=0 $
and by letting $\mu$ be the (usual) Liouville measure on the unit tangent
bundle $UB$. Thus we deduce that $\phi$ preserves this Liouville measure
$\mu$.  This proves \tref{quasi}. 

Note, that by construction  the Liouville measure $\mu$ is positive
on non-empty open subsets of $G$ and that the total mass of $\mu$
is proportional to the volume of $B$. The last one is finite if $M/\mathcal F$
is compact (\cite{BGP}).

\newpage

\section{Conjugate points} \label{secconj}
\subsection{Jacobi equation and Jacobi fields}  \label{jac}
We recall here some basic facts 
about the Jacobi equation, Jacobi fields and focal points. We refer
to \cite{jacobi} for extended explanations.

 Let $M$ be a Riemannian manifold, let $\gamma :I=[a,b]\to M$ be a
geodesic and let $\mathcal N$ be the normal bundle of $\gamma$.
Let $\Jac$ denote the space of all normal Jacobi fields along $\gamma$,
i.e. solutions of the equation $J'' + R(J) =0$, where $R$ denotes the
curvature endomorphism.  By $\omega$ we denote the canonical
symplectic form on $\Jac$, defined by 
$\omega (J_1,J_2)= \langle J_1 ' , J_2 \rangle +  \langle J_1  , J'_2 \rangle$.
For subspaces $W$ of $\Jac$ we denote by $W^{\perp}$ the orthogonal complement
with respect to $\omega$.  A subspace $W$ of $\Jac$ is called {\it isotropic}
if $W\subset W^{\perp}$ and it is called {\it Lagrangian} if $W=W^{\perp}$.
For an isotropic subspace $W$ and $t\in I$ we define the 
{\it $W$-focal index} of $t$ to be $f^W (t)= \dim (W) -\dim (W(t))$, where
$W(t) =\{ J(t)| J\in W\}$. The set of points with non-zero
focal index is discrete (\cite{jacobi})
and such points are called {\it $W$-focal}.
 The {\it $W$-index of $\gamma$}  is defined by 
$\ind _W (\gamma )= \Sigma _{t\in I} (f^W(t))$.  We have
the following semi-continuity property   (\cite{jacobi}):

\begin{lem} \label{cont}
Let $g_n$ be a sequence of Riemannian metrics that smoothly converges
to $g$. Let $\gamma _n :[a_n,b_n]\to M$ be a sequence of 
$g_n$-geodesics converging to $\gamma$. Let $W_n \subset \Jac (\gamma _n)$
be isotropic subspaces of normal Jacobi fields along $\gamma_n$ that
converge to an isotropic subspace $W\subset \Jac (\gamma )$. If 
$f^{W_n} (a_n) = f^W (a)$ and $f^{W_n} (b_n) = f^{W} (b)$ then
$\ind _{W_n} (\gamma _n) \leq \ind _{W_n} (\gamma)$ for all $n$ large
enough. If, in addition, $W_n$ are Lagrangians then 
$\ind _{W_n} (\gamma _n) = \ind _{W_n} (\gamma)$ for all $n$ large
enough.  
\end{lem}

 The following example is the main source of Lagrangians.
\begin{ex}
If $N$ is a submanifold of $M$ through $\gamma (a)$ orthogonal to $\gamma$,
then the space $\Lambda ^N$ of normal $N$-Jacobi fields is a Lagrangian. 
In this case the $\Lambda ^N$-focal index  of $a$ is equal to 
$\dim (M)-1 - \dim (N)$ and a point  $t\neq a$ 
is $\Lambda ^N$-focal if $\gamma (t)$  is a focal point of $N$ along 
$\gamma$ in the usual sense of Riemannian geometry.
In particular, the space $\Lambda = \Lambda ^{L(\gamma (a))}$
 of all $\mathcal F$-Jacobi fields along a 
horizontal geodesic $\gamma$ of a singular Riemannian foliation $\mathcal F$
is a Lagrangian. Thus the space $W$ of all vertical $\mathcal F$-Jacobi 
fields is isotropic.  
\end{ex}

 We recall now what we are going to use from Wilkings construction 
(\cite{Wilk})
of a transversal Jacobi equation. Let namely $W$ be an isotropic
space of Jacobi fields along $\gamma$.  Then there is a smooth Riemannian
vector bundle $\mathcal H$ with a Riemannian connection $'$ and
 with a Riemannian projection 
$P:\mathcal N  \to \mathcal H$, such that $\mathcal H(t) =\mathcal N(t) /W(t)$
for all $t$, that are not $W$-focal.  There is a smooth symmetric 
operator $R^{\mathcal H} :\mathcal H\to \mathcal H$ such 
that solutions of the  Jacobi equation $Y'' + R^{\mathcal H} (Y) =0$ are
precisely the projections (by the map $P$) of Jacobi fields 
$J\in W^{\perp}  \subset \Jac (\mathcal N)$ to $\mathcal H$.  Lagrangians
in $\Jac (\mathcal H)$ are precisely the projections of Lagrangians in
$\Jac (\mathcal N)$ that contain $W$.  Moreover (cf. \cite{jacobi}):

\begin{lem} \label{inddeco}
For each Lagrangian $\Lambda \subset \Jac (\mathcal N)$ that contains $W$,  we have
$\ind _W (\gamma)+ \ind _{\Lambda /W} (\gamma ) =\ind _{\Lambda } (\gamma )$.  
\end{lem} 

\begin{ex} \label{basicex}
  In the special case, where $\gamma$ is a horizontal geodesic with respect 
to a Riemannian submersion $f:M\to B$, let  $W$ be the space of 
$f$-vertical Jacobi fields, i.e., variational fields through variations of 
horizontal  lifts of $f(\gamma )$. Then $W(t)$ for each $t$ is the vertical
space of the submersion through $\gamma (t)$, $\mathcal H $ is canonically
identified with the normal bundle of the projected geodesic 
$\bar \gamma  =f(\gamma )$ in $B$ and the transversal operator 
$R^{\mathcal H}$ coincides with the curvature endomorphism in the base space
$B$.
\end{ex}

 \subsection{Vertical Jacobi fields} \label{finalsubsec}
Let $M$ be a Riemannian manifold and let $\mathcal F$
 be a singular Riemannian foliation on $M$.
Let $\gamma :[a,b] \to M$ be a horizontal geodesic. Then
the space
$\Lambda =\Lambda ^{L(\gamma (a))}$ 
of all normal $\mathcal F$-Jacobi fields along $\gamma$ is a
Lagrangian space of Jacobi fields. Note that the space 
$\Lambda$ depends not only on the maximal
geodesic containing $\gamma$ but also on the starting point
$\gamma (a)$.

 In the introduction, the space of $\mathcal F$-vertical Jacobi fields was
defined as the space of all Jacobi fields $J\in \Lambda$ with
$J(t) \in T_{\gamma (t)} L(\gamma (t))$ for all $t\in [a,b]$.
Recall now, that in Subsection \ref{vertjac} we have defined a space
$W^{\gamma}$ of Jacobi fields along $\gamma$, that are defined 
(independently of the starting point) as variational fields through
horizontal geodesics $\gamma _s$ with $\gamma _s (t) \in L (\gamma (t))$, for
all $t$. We have seen that 
$W^{\gamma } (t):= \{ J(t)| J\in W^{\gamma } \} $ coincides with
$ T_{\gamma (t)} L(\gamma (t))$, for all $t$. By definition 
 $W^{\gamma}\subset \Lambda$. Therefore,
$W^{\gamma}$ is precisely the space of all $\mathcal F$-vertical
Jacobi fields along $\gamma$.  In particular,  the latter
does not depend on the starting point, in contrast to $\Lambda$.

 The number $d(\gamma)$, defined in Subsection \ref{natdecomp},
is  the maximal dimension of $L(\gamma (t))$. We get 
$d(\gamma )= \dim W^{\gamma}$. The $W$-focal points along 
$\gamma$ are precisely the points $t_i$ with
$\dim L(\gamma (t_i)) < d(\gamma )$ and the $W$-focal index
of such points is  $d(\gamma ) -\dim L(\gamma (t_i))$.
In particular, for a  regular horizontal geodesic 
$\gamma$, its crossing number $c(\gamma )$ defined in the introduction
coincides with the {\it vertical index} $\ind _W (\gamma )$.

\begin{rem}
In the above terminology  it is possible to describe the space $W^{\perp}$
geometrically. Namely, it is possible to see that $W^{\perp}$ consists
of normal {\it horizontal} Jacobi fields, where we call a Jacobi field
$\mathcal F$-horizontal 
if it is the variational field of a variation of $\gamma$
through horizontal geodesics. This observation together with \lref{indeq} below
proves the equivalence between two a priori slightly different 
definitions of variational completeness used in \cite{Bott} and \cite{BottS}.
In our terminology  this equivalence reads as follows.  The singular Riemannian
foliation  $\mathcal F$ does not have horizontal conjugate points if and only 
if any $\mathcal F$-horizontal Jacobi field along any horizontal geodesic
that is tangent to the leaves at two points is tangent to the leaves at all
points. 
\end{rem}

\subsection{Horizontal conjugate points}
Let $\gamma$ be a curve  and let $T$ be a Riemannian bundle over
$\gamma$ with a Riemannian connection and a symmetric field of endomorphisms
$R:T\to T$ and the corresponding symplectic vector space $\Jac (T)$ of
Jacobi fields. (In this paper we are only interested in the cases 
$T=\mathcal N$ and $T=\mathcal H$).
Points $c<d$ in the interval of definition $I$ of $\gamma$ are called 
{\it conjugate} if there is a non-zero Jacobi fields $J\in \Jac$ with
$J(c)=J(d)=0$. Note that if $c<d$ are conjugate, then for 
each $\bar c \leq c$ there is some $\bar d \in [c,d]$ that is 
conjugate to $\bar c$ (\cite{jacobi}).

  Let now $M,\mathcal F,g$ be as always and let 
$\gamma :[a,b] \to M$ be  a horizontal geodesic.
Let $\Lambda =\Lambda ^{L(\gamma (a))}$ and 
$W=W^{\gamma}$ be defined as above.   Let $\mathcal H$ be the $W$-transversal
bundle as defined in Subsection \ref{jac}. The following result was 
independently obtained   in \cite{Nowak}:

\begin{lem} \label{indeq}
There are no horizontal conjugate points along $\gamma$ if and only if
$\ind _{\Lambda } (\gamma _0 ) =\ind _W (\gamma _0)$, where 
$\gamma _0$ denotes the subgeodesic $\gamma _0 : (a,b) \to M$ of $\gamma$.
  This condition is equivalent to
the statement that the point $a$ does not have conjugate points for the 
transversal Jacobi equation on $\mathcal H$.
\end{lem}
 
{\it Proof.}
Assume that  $\ind _{\Lambda } (\gamma ) =\ind _W (\gamma)$ and let some
$J \in  \Lambda$ and $a<t_0<b$
with $J(t_0) \in T_{\gamma (t_0)} L(\gamma (t_0))$ be given.
We find some $\tilde J \in W$, with $\tilde J(t_0)=J(t_0)$.
From   $\ind _{\Lambda } (\gamma ) =\ind _W (\gamma)$, we deduce
that $\tilde J -J \in W$ and therefore $J\in W$.
 The other implication is a direct consequence of the definition.

 To see the equivalence of $\ind _W (\gamma _0) =\ind _{\Lambda } (\gamma _0)$
to the absence of conjugate points in the quotient bundle $\mathcal H$,
we use \lref{inddeco} to see that  
$\ind _W (\gamma _0) =\ind _{\Lambda } (\gamma _0)$ is equivalent
to the absence of focal points of $\Lambda /W$ on the  open interval
$(a,b )$. But $\Lambda /W$ is by definition the Lagrangian  
$\tilde \Lambda ^a$  in $\Jac (\mathcal H)$ of all Jacobi fields $Y$
with $Y(a)=0$. Thus the statement that  $\ind _{\Lambda /W} (\gamma _0) =0$
is equivalent to the fact that $a$ does not have conjugate points with
respect to   the transversal Jacobi equation.
\qed

 Note that if the point $\gamma (a)$ is regular, then 
the condition $\ind _{\Lambda } (\gamma_0) =\ind _W (\gamma _0)$
is equivalent to $\ind _{\Lambda } (\gamma ) =
\ind _W (\gamma ) +\codim (\mathcal F ,M)-1$.

\subsection{Horizontal conjugate points in the  infinitesimally polar
case} Let $(M,g, \mathcal F)$ be as above and assume that
$\mathcal F$ is infinitesimally polar. Let $\gamma$ be a horizontal
geodesic in $M$. Then one can cover $\gamma$ by small
distinguished neighborhoods $O_i$. The restriction of $\mathcal F$ to each
$O_i$ has a Riemannian orbifold  $O_i /\mathcal F$ as quotient and
$\gamma \cap O_i$ is projected to an orbifold geodesic $\bar \gamma$
 in this quotient.
We get  a well defined development along the projection  $\bar \gamma$
of $\gamma$. 
We can consider it to be a smooth Riemannian manifold $B$
containing our geodesic $\gamma$ that we will denote by $\bar \gamma$,
if we consider it as part of $B$. 

\begin{rem} If $\mathcal F$ is infinitesimally
polar and closed and if 
$M$ is complete then $\gamma$ is projected to an orbifold-geodesic
$\bar \gamma$ in the Riemannian orbifold $M/\mathcal F$ and the manifold
$B$ considered above is the local development of $M/\mathcal F$ along 
$\bar \gamma$.
\end{rem}

 For regular geodesics, the following lemma is a direct consequence
of the  basic Example \ref{basicex}.

\begin{lem} \label{infpolconj}
There are no horizontal conjugate points along $\gamma$ if and only if
there are no conjugate points along $\bar \gamma$ in the local development 
$B$.
\end{lem}

{\it Proof.}
 Let $\gamma :[0,a] \to M$ be given.  
First let us assume that $b$ is horizontally conjugate to $0$ along $\gamma$.
Then there is a normal $L(\gamma (0))$-Jacobi field 
$J \in \Lambda \setminus W$ with
$J(b)=0$. Adding some element of the vertical space $W$ to $J$ we get
another $L(\gamma (0))$-Jacobi field $J_1$ with $J_1 (0) =0$ and $J_1 (b)\in 
T_{\gamma (b)}  L(\gamma (b))$.  Then $J_1$ is the variational field of
a variation $\gamma _s$ through horizontal geodesics  with $\gamma _0 =\gamma$
and $\gamma_s (0) =0$, for all $s$.  Then in each local
quotient $\gamma _s$ is projected to
a variation through  orbifold geodesics. Thus we obtain a lift $\bar \gamma _s$
to $B$ that gives us a variation of $\bar \gamma$ through geodesics in $B$
with $\bar \gamma _s (0)=0$. Moreover, the assumption 
$J_1 (b)=T_{\gamma (b)}L(\gamma (b))$ implies for the variational field 
$Y:= \frac d {ds} \bar \gamma _s$ that $Y(b)=0$.
If $b$ is non-conjugate to $a$ along $\bar \gamma$, then $Y$ is constant $0$.
But this is equivalent to $J\in W$.

  On the other hand, let us assume that there is some Jacobi field $Y$ along 
$\bar \gamma$ with $Y(0) =0$ and $Y(b)=0$. Find a variation of geodesics
$\bar \gamma _s$ corresponding to  this $Y$, with 
$\bar \gamma _s (0)=\gamma (0)$.  Since $T_{\gamma(0)} \mathcal F$
is polar we may choose a horizontal section $Z$ of  
$T_{\gamma(0)} \mathcal F$. Now we  find a unique smooth lift $\eta$
of the curve 
$\bar \eta (s): =\bar \gamma ' _s (0)$ to $Z$ with $\eta (0)=\gamma ' (0)$.
Then $\gamma _s (t)=\exp (t\eta (s))$ is a variation of $\gamma$ through 
horizontal geodesics, whose Jacobi field $J$ satisfies $J(0) =0$ and 
$J(b)\in T_{\gamma (b)} L(\gamma (b))$. Moreover, $J$ is not in $W$,
since $Y\neq 0$.
\qed

\begin{rem}
Since regular horizontal geodesics are dense in the space of all
geodesics and since the absence of conjugate points is an open condition
for Riemannian manifolds, we deduce from \lref{infpolconj} and
the proof of \pref{lastcor} below, that
a singular Riemannian  foliation $\mathcal F$ does not have horizontal
conjugate points if and only if any {\it regular} horizontal geodesic does not
have horizontal conjugate points.
\end{rem}

\subsection{Stability of absence of conjugate points}
Let $(M,g,\mathcal F)$ be as above. Let $\gamma :[a,b] \to M$  be a regular
horizontal geodesic. Let   
$\Lambda $ and $W$ be defined as in Subsection \ref{finalsubsec}. We assume
that $b$ is not $\Lambda$-focal. (The last condition can be achieved by
slightly increasing $b$).

 Let now $g_n$ be a sequence of Riemannian metrics on $M$ that smoothly 
converge  to $g$ and that are adapted to the singular Riemannian foliation
$\mathcal F$. Let $\gamma _n :[a,b] \to M$ be a sequence of $g_n$-horizontal
geodesics that converge to $\gamma$.  Let $\Lambda _n$ and $W_n$
be the spaces of $\mathcal F$-Jacobi fields and $\mathcal F$-vertical
Jacobi fields along $\gamma _n$ (with respect to  the metric $g_n$).

 Since we are in the regular
part of $\mathcal F$, the spaces $\Lambda _n$ converge to $\Lambda$
and $W_n$ converge to $W$. Moreover, $f^{W_n} =f^W(a)=f^{W_n} (b)=f^W(b)=0$
and $f^{\Lambda _n} (a)=f^{\Lambda } (a)= \dim (M)- 1 -\dim (\mathcal F )$
and $f^{\Lambda _n} (b)= f^{\Lambda} (b)=0$.
   From \lref{cont} we get
$\ind _{\Lambda } (\gamma ) -\ind _W (\gamma ) \leq \ind _{\Lambda _n} 
(\gamma _n) - \ind _{W_n} (\gamma _n)$, for all $n$ large enough.
Using \lref{indeq}   we derive:

\begin{lem}  \label{convergence}
In the above situation, assume that $\gamma _n$ has  no horizontal 
conjugate points. Then $\gamma$ has no horizontal conjugate points as well.
\end{lem}

\subsection{Conclusions}
Now we can finish the proofs of all results announced in the
introduction. We start with another characterization of infinitesimal
polarity.

\begin{prop} \label{lastcor}
 Let $\mathcal F$ be a singular Riemannian foliation on a Riemannian manifold
$M$ and  let $x$ be a point in $M$. Then $\mathcal F$ is infinitesimally
polar at $x$ if and only if there is a neighborhood $U$ of $x$,
such that the horizontal geodesics contained in $U$ do not have
horizontal conjugate points.
\end{prop}

{\it Proof.}
If $\mathcal F$ is infinitesimally polar at $x$
then we can find a small 
distinguished tubular neighborhood $O$ of $x$ such that $O/\mathcal F$
is a smooth Riemannian orbifold with bounded curvature. Then there 
is some $\epsilon >0$ such that each orbifold-geodesic 
of length $\leq \epsilon$  in $O/\mathcal F$ has no conjugate points.
Then taking $U\subset O$ to be an open  ball of radius $\epsilon$
around $x$, we deduce from   \lref{infpolconj} that no horizontal
geodesic in $U$ has horizontal conjugate points.

 On the other hand, let us assume that  $T_x \mathcal F$ is 
not polar. Then there is at least one regular horizontal geodesic
$\gamma$ in $T_x M$ with horizontal conjugate points, see \cite{LT}.
Using the convergence of the rescaled metrics on a small
tubular neighborhood $O$ of $x$ to the flat metric on $T_x M$
from Subsection \ref{infinit}, we deduce from \lref{convergence} 
that arbitrary small neighborhoods of $x$ contain regular geodesics 
$\gamma _n$ with horizontal conjugate points.
\qed

 Now we can finish the proof of  \tref{maingen}.  If there are no horizontal
conjugate points then $\mathcal F$ is infinitesimally polar by  \pref{lastcor}.
If, in addition, $\mathcal F$ is closed  then $B=M/\mathcal F$
is a complete Riemannian orbifold (Subsection \ref{rimorb}). From 
\lref{infpolconj}
we deduce that this Riemannian orbifold does not have conjugate points
if and only if in $M$ there are no horizontal conjugate points. 
Now the result follows from \lref{developp}.

\subsection{Continuity of the crossing counting function}
Now we are going to prove \tref{crosscont}.

{\it Proof.}
Recall that the crossing number $c(\gamma )$ of a regular geodesic is equal
to its vertical index $\ind _W (\gamma )$. The claim of \tref{crosscont} is local, i.e.,
$c$ is continuous if and only if each point $x\in M$ has a neighborhood
$U$ such that $c$ is continuous for the restricted singular Riemannian
foliation $(U,\mathcal F)$.

  If $\mathcal F$ is infinitesimally polar at $x$, we may choose
$U$ as in \pref{lastcor}. Then for each regular geodesic 
$\gamma$ in $U$ we get 
$\ind _W(\gamma )=\ind _{\Lambda} (\gamma )- (\codim (\mathcal F, M) -1 )$. 
The result
now follows from the continuity of indices for Lagrangians, \lref{cont}.

Now let us assume that $c$ is continuous and that 
$T_x \mathcal F$ is not polar.
The proof of \cref{lastcor} shows that there are regular geodesics
$\gamma _n :[0,\epsilon _n] \to M$, with $\epsilon _n \to 0$ and
$\gamma _n (0) \to x$  that have horizontal conjugate points and such that
the starting directions $v_n $ of $\gamma _n$ converge to a regular
direction $v$ in $T_x M$. Let $\gamma=\gamma _v$  be the geodesic in $M$
 in the direction 
$v$.  For sufficiently small
$\epsilon=\epsilon (\gamma )$, 
the geodesic $\gamma :[-\epsilon ,\epsilon ] \to M$ has
no horizontal conjugate points.  But the extended geodesics 
$\gamma _n :[-\epsilon, \epsilon] \to M$ still have horizontal conjugate 
points. 
Thus, for $n$ large enough,  
$c(\gamma _n)= \ind _{W_n } (\gamma _n) < \ind _{\Lambda _n}  (\gamma _n)
-(\codim (\mathcal F ,M) -1)$ and $c(\gamma))= \ind _W (\gamma )=  
\ind _{\Lambda } (\gamma )- (\codim (\mathcal F ,M) -1)$.
Since $\gamma _n$ converges to $\gamma$,
this contradicts  the continuity of $c$.
\qed

\bibliographystyle{alpha}
\bibliography{expl}

\end{document}